\documentclass[11pt,reqno]{amsart}
\usepackage{latexsym,amsmath,amssymb,amsfonts,amscd,amsthm}
\usepackage{color}
\makeatletter
\theoremstyle{plain}
\newtheorem{thm}{Theorem}
\theoremstyle{plain}

\numberwithin{equation}{section}
\numberwithin{figure}{section} 
\theoremstyle{plain}
\newtheorem{cor}{Corollary}
\theoremstyle{plain}
\newtheorem{lem}{Lemma}
\theoremstyle{plain}

\theoremstyle{plain}

\theoremstyle{plain}

\theoremstyle{plain}

\theoremstyle{plain}

\theoremstyle{plain}
\newtheorem{definition}{Definition}

\begin{document}

\begin{center}

\textbf{\large{Third order differential subordination and superordination results for analytic  \\[1mm] functions
 involving the Srivastava-Attiya operator}}\\[4mm]

\textbf{H. M. Srivastava$^{1,2,\ast}$, 
A. Prajapati$^{3}$ and P. Gochhayat$^{4}$}\\[2mm]
 
$^{1}$Department of Mathematics and Statistics,
University of Victoria, \\

Victoria, British Columbia V8W 3R4, Canada\\[1mm]

$^{2}$Department of Medical Research, China Medical University Hospital,\\

China Medical University, Taichung 40402, Taiwan, Republic of China\\[1mm]

\textbf{E-Mail: harimsri@math.uvic.ca}\\[1mm]

\textbf{$^{\ast}$Corresponding Author}\\[2mm]

$^{3,4}$Department of Mathematics , \\

Sambalpur University, Jyoti vihar 768019, Burla ,\\
Sambalpur, Odisha, India \\[1mm]
\textbf {$^{3}$ E-Mail: anujaprajapati49@gmail.com}\\[1mm]
\textbf {$^{4}$ E-Mail: pgochhayat@gmail.com}\\[2mm]

\textbf{Abstract}
\end{center}
\begin{quotation}
In this article, by making use of the linear operator introduced and studied by Srivastava and Attiya \cite{srivastava1}, suitable classes of admissible functions are investigated and the dual properties of the third-order differential subordinations are presented. As a consequence, various sandwich-type theorems are established for a class of univalent analytic functions involving the celebrated  Srivastava-Attiya transform. Relevant connections of the new results are pointed out.\\
\end{quotation}



\section{\bf Introduction, Definitions and Preliminaries}
Let $\mathcal{H}$  be  the class  of  functions  analytic  in the
\textit{open} unit disk $$
\mathbb {U}:=\{z:z\in\mathbb{C}~\text{and}~|z|< 1\}.$$
Also let $$\mathcal{H}[a,n]\quad (n \in \mathbb{N}:=\{1,2,3,\cdots\},a \in \mathbb{C})$$  be the subclass of the analytic function class $\mathcal{H}$ consisting of functions of the form
$$f(z)=a+a_{n}z^{n}+a_{n+1}z^{n+1}+\cdots, \quad(z \in \mathbb{U}).$$ Let $\mathcal{A}(\subset \mathcal{H})$ be the class of  functions which are analytic in $\mathbb U$  and have the \emph{normalized} Taylor-Maclaurin  series of the form:
\begin{equation}\label{z1}
\mathit{f}(z)=z+\sum_{n=2}^{\infty}a_{n}z^{n},\qquad(z \in \mathbb{U}).
\end{equation}
Suppose that $f$ and $g$ are in $\mathcal H$. We say that $f$ is \textit{subordinate} to $g$, (or $g$ is
\textit{superordinate} to $f$), written as $$f \prec g ~\text{
	in}~\mathbb U~\text{or}~ f(z) \prec g(z) \quad (z \in
\mathbb U),$$ if there exists a function $\omega\in \mathcal H$,
satisfying the conditions of the Schwarz lemma, namely
$$~\omega(0)=0 ~\text{and}~ |\omega(z)|< 1$$ such that
$$f(z)=g(\omega(z))\quad(z \in \mathbb U).$$ It follows that
$$
f(z) \prec g(z)\;(z \in \mathbb U) \Longrightarrow f(0)=g(0) ~
\text{and} ~ f(\mathbb U) \subset g(\mathbb U).
$$ In particular, if $g$ is \textit{univalent} in $\mathbb U$, then
the reverse implication also holds (see, for details, \cite{miller3}).

The concept of differential subordination is a generalization of various inequalities involving complex variables. We recall here some more definitions and terminologies from the theory of differential subordination and  superordination.
\begin{definition}\label{d3a}\cite{antonion}
	{\rm Let $\psi:\mathbb{C}^{4} \times \mathbb{U}\longrightarrow \mathbb{C}$ and suppose that the function $h(z)$ is univalent in $\mathbb{U}.$ If the function $p(z)$ is analytic in $\mathbb{U}$ and satisfies the following third-order differential subordination
	\begin{equation}\label{p1}
	\psi(p(z),zp^{\prime}(z),z^{2}p^{\prime\prime}(z),z^{3}
	p^{\prime\prime\prime}(z);z)\prec h(z),
	\end{equation}
	then $p(z)$ is called a \emph{solution} of the differential subordination \eqref{p1}. Furthermore, a given univalent function $q(z)$ is called a \emph{dominant of the solutions of the differential subordination \eqref{p1}, or, more simply, a dominant} if $p(z)\prec q(z)$ for all $p(z)$ satisfying (\ref{p1}). A dominant $\tilde{q}(z)$ 
	that satisfies $\tilde{q}(z)\prec q(z)$ 
	for all dominants $q(z)$ of (\ref{p1}) is said to be the \emph{best dominant}.}.
\end{definition}
 
\begin{definition}\label{d5}\cite{tang}
	{\rm Let $\psi:\mathbb{C}^{4} \times \mathbb{U}\longrightarrow \mathbb{C}$ and the function $h(z)$ be univalent in $\mathbb{U}.$ If the function $p(z)$ and
	$$\psi(p(z),zp^{\prime}(z),z^{2}p^{\prime\prime}(z),z^{3}
	p^{\prime\prime\prime}(z);z)$$
	are univalent in $\mathbb{U}$ and satisfies the following third-order differential superordination
	\begin{equation}\label{p2}
	h(z) \prec \psi(p(z),zp^{\prime}(z),z^{2}p^{\prime\prime}(z),z^{3}
	p^{\prime\prime\prime}(z);z),
	\end{equation}
	then $p(z)$ is called a \emph{solution} of the differential superordination. An analytic function $q(z)$ is called a \emph{subordinant} of the solutions of the differential superordination, or more simply a subordinant, if $ q(z) \prec p(z)$ for all $p(z)$ satisfying (\ref{p2})}.
\end{definition}
A univalent subordinant $\tilde{q}(z)$ that satisfies $q(z)\prec \tilde{q}(z)$ 
for all subordinant $q(z)$ of (\ref{p2}) is said to be the \emph{best subordinant}. We note that both the best dominant and best subordinant are unique up to rotation of $\mathbb U$. The well known monograph
of Miller and Mocanu \cite{miller3} and the more recent book of Bulboac\u{a} \cite{bulboaca} provide detailed
expositions on the theory of differential subordination and superordination.
With a view to define the Srivastava-Attiya  transform we recall here the generalized Hurwitz-Lerch Zeta function, which is defined in \cite{srivastava2} by the following series:
\begin{gather}\label{g3}
\Upsilon(z,\mu,b)=\sum_{n=0}^{\infty}\frac{z^{n}}{(b+n)^{\mu}},\\~~\left(b \in \mathbb{C}\setminus \mathbb{Z}_{0}^{-}; \mu \in \mathbb{C} ~~\text{when}~~ z \in \mathbb{U}; \Re (\mu)>1 ~~\text{when}~~ z \in \partial \mathbb{U}\right)\nonumber.
\end{gather}
Special cases of the function $\Upsilon (z,\mu,b)$ include for example, the Riemann Zeta function $\zeta(\mu)=\Upsilon(1,\mu,1)$; the Hurwitz Zeta function $\zeta(\mu,b)=\Upsilon(1,\mu,b),$ the Lerch Zeta function $l_{\mu}{\zeta}=\Upsilon(\exp 2\pi i\zeta,\mu,1)\\ (\zeta \in \mathbb{R}, \Re(\mu)>1)$, the Poly logarithm function $L_{i\mu}=z \Upsilon(z,\mu,1)$ and so on. For further details see \cite{srivastava3} and the references therein.
Srivastava and Attiya  \cite{srivastava1} considered the following normalized function:
\begin{equation}\label{g2}
R_{\mu,b}(z)=(1+b)^{\mu}[\Upsilon(z,\mu,b)-b^{-\mu}]=z+\sum_{n=2}^{\infty}\left(\frac{b+1}{b+n}\right)^{\mu}z^{n}, \qquad (z \in \mathbb{U}),
\end{equation}
 and by making use of $R_{\mu,b}(z)$, they have introduced the  linear operator $J_{\mu,b}:\mathcal{A} \rightarrow \mathcal{A}$ which is defined in terms of convolution as follows:
\begin{equation}\label{g4}
J_{\mu,b}f(z)=R_{\mu,b}(z)*f(z)=z+\sum_{n=2}^{\infty}\left(\frac{b+1}{b+n}\right)^{\mu} a_{n}z^{n},\qquad (z \in \mathbb{U}).
\end{equation}
The operator $J_{\mu,b}f(z)$ is now popularly known in the literature as the Srivastava-Attiya operator. Various applications  of $J_{\mu,b}f(z)$  are found in \cite{cho3,pg1,pg56,pg,pg2,wang} and the references therein. 
 From  (\ref{g4}), it is clear that
\begin{equation}\label{g5}
zJ_{\mu+1,b}^{\prime}f(z)=(b+1)J_{\mu,b}f(z)-bJ_{\mu+1,b}f(z).
\end{equation} 
Suitable choices of parameters, the above defined operator unifies various other linear operators which are introduced earlier. For examples
\begin{enumerate}
\item $ J_{0,b}f(z)=f(z)$,
\item $J_{1,0}f(z)= \int _{0}^{z}\frac{f(t)}{t}dt :=\mathfrak{A} f(z),$
\item $J_{1,\eta}f(z)=\frac{1+\eta}{z^{\eta}} \int_{0}^{z}t^{\eta-1}f(t)dt :=\mathfrak{I}_{\eta}f(z),\qquad (\eta>-1), $
\item $J_{\sigma,1}f(z)=z+\sum _{n=2}^{\infty}\left(\frac{2}{n+1}\right)^{\sigma}a_{n}z^{n}:= I^{\sigma}f(z)\qquad (\sigma>0).$
\end{enumerate} 
Where $\mathfrak{A}(f)$ and  $\mathfrak{I}_{\eta}$ are the integral operators introduced by Alexander  and Bernardi, respectively, and $I^{\sigma}(f)$ is the Jung-Kim-Srivastava integral operator  closely related to multiplier transformation studied by Flett. For more detail unifications  we refer \cite{pg2}.

\begin{definition}\label{d4}\cite{antonion}
	{\rm Let $\mathbb{Q}$ be the set of all  functions $q$ that are analytic and univalent on  $\overline{\mathbb{U}}\setminus E(q),$ where $E(q)=\{\xi:\xi \in \partial \mathbb{U}:\lim_{z\rightarrow \xi}q(z)= \infty\},$ and are such that $\min \mid q^{\prime}(\xi) \mid=\rho >0$ for $\xi \in \partial \mathbb{U} \setminus E(q).$
	Further, let the subclass of $\mathbb{Q}$ for which $q(0)=a$ be denoted by $\mathbb{Q}(a),\mathbb{Q}(0)=\mathbb{Q}_{0}$ and $\mathbb{Q}(1)=\mathbb{Q}_{1}.$}
\end{definition}
The subordination methodology is applied to an appropriate class of admissible functions. The following class of admissible functions is given by Antonino and Miller.
\begin{definition}\label{d2a}\cite{antonion}
	{\rm Let $\Omega$ be a set in $\mathbb{C}$ and $q \in \mathbb{Q}$  and $n \in \mathbb{N} \setminus \{1\}.$ The class of admissible functions $\Psi_{n}[\Omega,q]$ consists of those functions $\psi:\mathbb{C}^{4} \times \mathbb{U}\longrightarrow \mathbb{C}$ achieving the following admissibility conditions:
	$$\psi(r,s,t,u;z) \not\in \Omega$$
	whenever
	$$r=q(\zeta),s=k\zeta q^{\prime}(\zeta), \Re\left(\frac{t}{s}+1\right)\geq k \Re\left(\frac{\zeta q^{\prime\prime}(\zeta)}{q^{\prime}(\zeta)}+1\right),$$
	and $$\Re\left(\frac{u}{s}\right)\geq k^{2} \Re\left(\frac{\zeta^{2} q^{\prime\prime\prime}(\zeta)}{q^{\prime}(\zeta)}\right),$$
	where $z \in \mathbb{U},\zeta \in \partial \mathbb{U} \setminus E(q),$ and $k\geq n.$}
\end{definition}
The next lemma is the foundation result in the theory of third-order differential subordination.
\begin{lem}\label{t4}\cite{antonion}
	Let $p \in \mathcal{H}[a,n]$ with $n\geq 2,$ and $q \in \mathbb{Q}(a)$ achieving the following conditions:
	$$\Re\left(\frac{\zeta q^{\prime\prime}(\zeta)}{q^{\prime}(\zeta)}\right)\geq 0,\quad \left|\frac{zq^{\prime}(z)}{q^{\prime}(\zeta)}\right|\leq k,$$ where $z \in \mathbb{U},\zeta \in \partial \mathbb{U} \setminus E(q),$ and $k\geq n.$ If $\Omega$ is a set in $\mathbb{C},\psi \in \Psi_{n}[\Omega,q]$ and $$\psi \left(p(z),zp^{\prime}(z),z^{2}p^{\prime\prime}(z),z^{3}
	p^{\prime\prime\prime}(z);z\right) \subset \Omega,$$
	then $$p(z) \prec q(z)\qquad (z \in \mathbb{U}).$$
\end{lem}
\begin{definition}\label{d6}\cite{tang}
	{\rm Let $\Omega$ be a set in $\mathbb{C}$ and $q \in \mathcal{H}[a,n]$  and $ q^{\prime}(z)\neq 0.$ The class of admissible functions $\Psi_{n}^{\prime}[\Omega,q]$ consists of those functions $\psi:\mathbb{C}^{4} \times \mathbb{\overline{U}}\longrightarrow \mathbb{C}$ that satisfy the following admissibility conditions:
	$$\psi(r,s,t,u;\zeta) \in \Omega$$
	whenever
	$$r=q(z),s=\frac{zq^{\prime}(z)}{m}, \Re\left(\frac{t}{s}+1\right)\leq \frac{1}{m} \Re\left(\frac{z q^{\prime\prime}(z)}{q^{\prime}(z)}+1\right),$$
	and $$\Re\left(\frac{u}{s}\right)\leq \frac{1}{m^{2}} \Re\left(\frac{z^{2} q^{\prime\prime\prime}(z)}{q^{\prime}(z)}\right),$$
	where $z \in \mathbb{U},\zeta \in \partial \mathbb{U},$ and $m\geq n \geq 2.$}
\end{definition} 
\begin{lem}\label{t5}\cite{tang}
	Let $p \in \mathcal{H}[a,n]$ with $\psi \in \Psi_{n}^{\prime}[\Omega,q].$ If
	$$\psi(p(z),zp^{\prime}(z),z^{2}p^{\prime\prime}(z),z^{3}
	p^{\prime\prime\prime}(z);z)$$ is univalent in $\mathbb{U}$ and $p \in \mathbb{Q}(a)$
	satisfying the following conditions:
	$$\Re\left(\frac{z q^{\prime\prime}(z)}{q^{\prime}(z)}\right)\geq 0,\quad\left|\frac{zp^{\prime}(z)}{q^{\prime}(z)}\right|\leq m,$$ where $z \in \mathbb{U},\zeta \in \partial \mathbb{U},$ and $m\geq n \geq 2,$ then $$\Omega \subset \{\psi \left(p(z),zp^{\prime}(z),z^{2}p^{\prime\prime}(z),z^{3}p^{\prime\prime\prime}(z);z\right):z \in \mathbb{U}\},$$
implies that $$q(z) \prec p(z)\qquad  (z \in \mathbb{U}).$$
\end{lem}
Though the notion of third order differential subordination have originally found in the work of Ponnusamy and Juneja \cite{ponnusamy}.   The recent work due to Tang et al. \cite{tang,tang2} on third order differential  subordination attracted to many researchers in this field. For example see \cite{farzana,ibrahim,jeyaraman, raducanu,tang,tang1,tang2,ponnusamy,miller4}. In the present paper we  considered suitable classes of admissible functions associated with Srivastava-Attiya operator and   obtained sufficient conditions on the normalized analytic function $f$ such that \textit{Sandwich-type subordination} of the following form holds:
 $$h_{1}(z)\prec \Theta(f) \prec q_{2}(z),\qquad (z \in \mathbb{U}),$$
 where $q_{1},q_{2}$ are univalent in $\mathbb{U}$ and $\Theta$ is a suitable operator.
\section{\bf  Results Related to Third Order Subordination}
In this section, start with given set $\Omega$ and given function $q$ and we determine a set of admissible operators $\psi$ so that (\ref{p1}) holds true. Thus, the following new class of admissible function is introduced which will required to prove the main third-order differential subordination theorems for the operator $J_{\mu,b}f(z)$  defined by (\ref{g2}).
\begin{definition}\label{d1}
{\rm Let $\Omega$ be a set in $\mathbb{C}$ and $q \in \mathbb{Q}_{0}\bigcap \mathcal{H}_{0}.$ The class of admissible function $\Phi_{J}[\Omega,q]$ consists of those functions $\phi:\mathbb{C}^{4}\times \mathbb{U}\longrightarrow\mathbb{C}$ that satisfy the following admissibility conditions:
\begin{equation*}
\phi(\alpha,\beta,\gamma,\delta;z) \not\in \Omega
\end{equation*}
whenever
\begin{equation*}
\alpha=q(\zeta), \beta=\frac{k\zeta q^{\prime}(\zeta)+bq(\zeta)}{b+1},$$ $$\Re\left(\frac{\gamma(b+1)^{2}-b^{2}\alpha}{(\beta(b+1)-b\alpha)}-2b\right)\geq k \Re\left(\frac{\zeta q^{\prime\prime}(\zeta)}{q^{\prime}(\zeta)}+1\right),
\end{equation*}
and 
\begin{equation*}
\Re\left(\frac{\delta(b+1)^{3}-\gamma(b+1)^{2}(3b+3)+b^{2}\alpha(3+2b)}{(b(\beta-\alpha)+\beta)}+3b^{2}+6b+2\right)
\geq k^{2} \Re\left(\frac{\zeta^{2} q^{\prime\prime\prime}(\zeta)}{q^{\prime}(\zeta)}\right),
\end{equation*}
where $z \in \mathbb{U},\zeta \in \partial \mathbb{U} \setminus E(q),$ and $k\geq 2.$}
\end{definition}
\begin{thm}\label{t6}
Let $\phi \in \Phi_{J}[\Omega,q].$ If the function $f \in \mathcal{A}$ and $q \in \mathbb{Q}_{0}$ satisfy the following conditions:
\begin{equation}\label{p3}
\Re\left(\frac{\zeta q^{\prime\prime}(\zeta)}{q^{\prime}(\zeta)}\right)\geq 0,\quad \left|\frac{J_{\mu,b}f(z)}{q^{\prime}(\zeta)}\right|\leq k,
\end{equation}
and 
\begin{equation}\label{p4}
\{\phi(J_{\mu+1,b}f(z),J_{\mu,b}f(z),J_{\mu-1,b}f(z),J_{\mu-2,b}f(z)
;z):z \in \mathbb{U}\}\subset \Omega,
\end{equation}
then 
\begin{equation*}
J_{\mu+1,b}f(z) \prec q(z)\qquad(z \in \mathbb{U}).
\end{equation*}
\begin{proof}
Define the analytic function $p(z)$ in $\mathbb{U}$ by
\begin{equation}\label{p5}
p(z)=J_{\mu+1,b}f(z).
\end{equation}
From equation (\ref{g5}) and (\ref{p5}), we have
\begin{equation}\label{p6}
J_{\mu,b}f(z)=\frac{zp^{\prime}(z)+bp(z)}{b+1}.
\end{equation}
By similar argument, yields
\begin{equation}\label{p7}
J_{\mu-1,b}f(z)=\frac{z^{2}p^{\prime\prime}(z)+(2b+1)zp^{\prime}(z)+b^{2}p(z)}{(b+1)^{2}}
\end{equation}
and
\begin{equation}\label{p8}
J_{\mu-2,b}f(z)=\frac{z^{3}p^{\prime\prime\prime}(z)+(3b+3)z^{2}p^{\prime\prime}(z)+(3b^{2}+3b+1)zp^{\prime}(z)+b^{3}p(z)}{(b+1)^{3}}.
\end{equation}
Define the transformation from $\mathbb{C}^{4}$ to $\mathbb{C}$ by
\begin{equation*}
\alpha(r,s,t,u)=r,\qquad \beta(r,s,t,u)=\frac{s+br}{b+1},
\end{equation*}
\begin{equation}\label{p9}
\gamma(r,s,t,u)=\frac{t+(2b+1)s+b^{2}r}{(b+1)^{2}}
\end{equation}
and
\begin{equation}\label{p10}
\delta(r,s,t,u)=\frac{u+(3b+3)t+(3b^{2}+3b+1)s+b^{3}r}{(b+1)^{3}}.
\end{equation}
Let
\begin{multline}\label{p11}
\psi(r,s,t,u)=\phi(\alpha,\beta,\gamma,\delta;z)=
\phi \bigg(r,\frac{s+br}{b+1},\frac{t+(2b+1)s+b^{2}r}{(b+1)^{2}},\\  \frac{u+(3b+3)t+(3b^{2}+3b+1)s+b^{3}r}{(b+1)^{3}};z\bigg).
\end{multline}
The proof will make use of  Lemma \ref{t4}. Using equations (\ref{p5}) to (\ref{p8}),  and from (\ref{p11}), we have
\begin{equation}\label{p12}
\psi \left(p(z),zp^{\prime}(z),z^{2}p^{\prime\prime}(z),z^{3}p^{\prime\prime\prime}(z);z\right)=\phi \left(J_{\mu+1,b}f(z),J_{\mu,b}f(z),J_{\mu-1,b}f(z),J_{\mu-2,b}f(z)
;z \right).
\end{equation} 
Hence, (\ref{p4}) becomes
\begin{equation*}
\psi \left(p(z),zp^{\prime}(z),z^{2}p^{\prime\prime}(z),z^{3}p^{\prime\prime\prime}(z);z\right) \in \Omega.
\end{equation*}
Note that
\begin{equation*}
\frac{t}{s}+1=\frac{\gamma(b+1)^{2}-b^{2}\alpha}{(\beta(b+1)-b\alpha)}-2b
\end{equation*}
and
\begin{equation*}
\frac{u}{s}=\frac{\delta(b+1)^{3}-\gamma(b+1)^{2}(3b+3)+b^{2}\alpha(3+2b)}{(b(\beta-\alpha)+\beta)}.
\end{equation*}
Thus, the admissibility condition for $\phi \in \Phi_{J}[\Omega,q]$ in Definition \ref{d1} is equivalent to the admissibility condition for $\psi \in \Psi_{2}[\Omega,q]$ as given in Definition \ref {d2a} with $n=2.$ Therefore, by using (\ref{p3}) and Lemma \ref{t4}, we have
\begin{equation*}
J_{\mu+1,b}f(z)\prec q(z).
\end{equation*}
This completes the proof of theorem.
\end{proof}
\end{thm}
The next result is an extension of Theorem \ref{t6}  to the case where the behavior of $q(z)$ on $\partial \mathbb{U}$ is not known.
\begin{cor}\label{c1}
Let $\Omega \subset \mathbb{C}$ and let the function $q$ be univalent in $\mathbb{U}$ with $q(0)=0.$ Let $\phi \in \Phi_{J}[\Omega,q_{\rho}]$ for some $\rho \in (0,1),$ where $q_{\rho}(z)=q(\rho z).$ If the function $f \in \mathcal{A}$ and $q_{\rho}$ satisfy the following conditions
\begin{equation*}
\Re\left(\frac{\zeta q_{\rho}^{\prime\prime}(\zeta)}{q_{\rho}^{\prime}(\zeta)}\right)\geq 0,\quad\left|\frac{J_{\mu,b}f(z)}{q_{\rho}^{\prime}(\zeta)}\right|\leq k\qquad(z \in \mathbb{U}, k \geq 2, \zeta \in \partial \mathbb{U}\setminus E(q_{\rho}))
\end{equation*}
and
\begin{equation*}
\phi(J_{\mu+1,b}f(z),J_{\mu,b}f(z),J_{\mu-1,b}f(z),J_{\mu-2,b}f(z);z) \in \Omega,
\end{equation*}
then
\begin{equation*}
J_{\mu+1,b}f(z)\prec q(z)\qquad(z \in \mathbb{U}).
\end{equation*}
\begin{proof}
From Theorem \ref{t6}, then $J_{\mu+1,b}f(z)\prec q_{\rho}(z).$ The result asserted by Corollary \ref{c1} is now deduced from the following subordination property $q_{\rho}(z)\prec q(z)\qquad(z \in \mathbb{U}).$ 
\end{proof}
\end{cor}
If $\Omega \neq \mathbb{C}$ is a simply connected domain, then $\Omega=h(\mathbb{U})$ for some conformal mapping $h(z)$ of $\mathbb{U}$ onto $\Omega.$ In this case, the class $\Phi_{J}[h(\mathbb{U}),q]$ is written  as $\Phi_{J}[h,q].$ This follows immediate consequence of Theorem  \ref{t6}.
\begin{thm}\label{t7}
Let $\phi \in \Phi_{J}[h,q].$ If the function $f \in \mathcal{A}$ and $q \in \mathbb{Q}_{0}$ satisfy the following conditions:
\begin{equation}\label{p13}
\Re\left(\frac{\zeta q^{\prime\prime}(\zeta)}{q^{\prime}(\zeta)}\right)\geq 0, \quad \left|\frac{J_{\mu,b}f(z)}{q^{\prime}(\zeta)}\right|\leq k,
\end{equation}
and 
\begin{equation}\label{p14}
\phi(J_{\mu+1,b}f(z),J_{\mu,b}f(z),J_{\mu-1,b}f(z),J_{\mu-2,b}f(z);z)\prec h(z),
\end{equation}
then 
\begin{equation*}
J_{\mu+1,b}f(z) \prec q(z)\qquad (z \in \mathbb{U}).
\end{equation*}
\end{thm}
The next result is an immediate consequence of Corollary \ref{c1}.
\begin{cor}\label{c2}
Let $\Omega \subset \mathbb{C}$ and let the function $q$ be univalent in $\mathbb{U}$ with $q(0)=0.$ Let $\phi \in \Phi_{J}[h,q_{\rho}]$ for some $\rho \in (0,1),$ where $q_{\rho}(z)=q(\rho z).$ If the function $f \in \mathcal{A}$ and $q_{\rho}$ satisfy the following conditions
\begin{equation*}
\Re\left(\frac{\zeta q_{\rho}^{\prime\prime}(\zeta)}{q_{\rho}^{\prime}(\zeta)}\right)\geq 0,\quad \left|\frac{J_{\mu,b}f(z)}{q_{\rho}^{\prime}(\zeta)}\right|\leq k \qquad (z \in \mathbb{U}, k \geq 2, \zeta \in \partial \mathbb{U}\setminus E(q_{\rho})),
\end{equation*}
and
\begin{equation*}
\phi(J_{\mu+1,b}f(z),J_{\mu,b}f(z),J_{\mu-1,b}f(z),J_{\mu-2,b}f(z);z) \prec h(z),
\end{equation*}
then
\begin{equation*}
J_{\mu+1,b}f(z)\prec q(z)\qquad (z \in \mathbb{U}).
\end{equation*}
\end{cor}
The following result  yields the best dominant of the differential subordination (\ref{p14}).
\begin{thm}\label{t8}
Let the function $h$ be univalent in $\mathbb{U}$ and let $\phi : \mathbb{C}^{4}\times \mathbb{U}\longrightarrow \mathbb{C}$ and $\psi$ be given by (\ref{p11}). Suppose that the differential equation
\begin{equation}\label{p15}
\psi(q(z),zq^{\prime}(z),z^{2}q^{\prime\prime}(z),z^{3}
q^{\prime\prime\prime}(z);z)=h(z),
\end{equation}
has a solution $q(z)$ with $q(0)=0,$ which satisfy condition (\ref{p3}). If the function $f \in \mathcal{A}$ satisfies condition (\ref{p14}) and 
\begin{equation*}
\phi(J_{\mu+1,b}f(z),J_{\mu,b}f(z),J_{\mu-1,b}f(z),J_{\mu-2,b}f(z);z)
\end{equation*}
is analytic in $\mathbb{U},$ then
\begin{equation*}
J_{\mu+1,b}f(z)\prec q(z)
\end{equation*}
and $q(z)$ is the best dominant.
\begin{proof}
From Theorem \ref{t6}, we have $q$ is a dominant of (\ref{p14}). Since $q$ satisfies (\ref{p15}), it is also a solution of (\ref{p14}) and therefore $q$ will be dominated by all dominants. Hence $q$ is the best dominant.
\end{proof}
\end{thm}
In view of Definition \ref{d1}, and in the special case $q(z)=Mz,\quad M>0,$ the class of admissible functions $\Phi_{J}[\Omega,q],$ denoted by  $\Phi_{J}[\Omega,M],$ is expressed as follows.
\begin{definition}\label{d2}
{\rm Let $\Omega$ be a set in $\mathbb{C}$ and $M>0.$ The class of admissible function $\Phi_{J}[\Omega,M]$ consists of those functions $\phi:\mathbb{C}^{4}\times \mathbb{U}\longrightarrow\mathbb{C}$ such that
\begin{multline}\label{p16}
\phi \bigg(Me^{i\theta},\frac{(k+b)Me^{i\theta}}{b.+1},\frac{L+[(2b+1)k+b^{2}]Me^{i\theta}}{(b+1)^{2}},
\frac{N+(3b+3)L+[(3b^{2}+3b+1)k+b^{3}]Me^{i\theta}}{(b+1)^{3}};z\bigg)\notin \Omega.
\end{multline}
whenever $z \in \mathbb{U},\Re(Le^{-i\theta})\geq (k-1)kM,$ and $\Re (Ne^{-i\theta})\geq 0$ for all $\theta \in \mathbb{R}$ and $k \geq 2.$ }
\end{definition}
\begin{cor}\label{c3}
Let $\phi \in \Phi_{J}[\Omega,M].$ If the function $f \in \mathcal{A}$ satisfies 
\begin{equation*}
\left|J_{\mu,b}f(z)\right|\leq kM\qquad( z \in \mathbb{U},k\geq 2 ; M>0),
\end{equation*}
and 
\begin{equation*}
\phi(J_{\mu+1,b}f(z),J_{\mu,b}f(z),J_{\mu-1,b}f(z),J_{\mu-2,b}f(z);z)\in \Omega,
\end{equation*}
then
\begin{equation*}
\left|J_{\mu+1,b}f(z)\right|<M.
\end{equation*}
\end{cor}
 In this special case $\Omega=q(\mathbb{U})=\{w:|w|<M\},$ the class $\Phi_{J}[\Omega,M]$ is simply denoted by $\Phi_{J}[M].$ Corollary \ref{c3} can now be written in the following form.
 \begin{cor}\label{c4}
 Let $\phi \in \Phi_{J}[M].$ If the function $f \in \mathcal{A}$ satisfies 
\begin{equation*}
\left|J_{\mu,b}f(z)\right|\leq kM\qquad(k\geq 2; z \in \mathbb{U}, M>0),
\end{equation*}
and 
\begin{equation*}
\left|J_{\mu+1,b}f(z),J_{\mu,b}f(z),J_{\mu-1,b}f(z),J_{\mu-2,b}f(z);z \right|<M,
\end{equation*}
then
\begin{equation*}
\left|J_{\mu+1,b}f(z)\right|< M.
\end{equation*}
\end{cor}
\begin{cor}\label{s3}
Let $k\geq 2,\quad 0 \neq \mu \in \mathbb{C}$ and $M>0.$ If the function $f \in \mathcal{A}$ satisfies
\begin{equation*}
|J_{\mu,b}f(z)|\leq k M,
\end{equation*}
and
\begin{equation*}
\left|J_{\mu,b}f(z)-J_{\mu+1,b}f(z)\right|<\frac{M}{|b+1|},
\end{equation*}
then
\begin{equation*}
|J_{\mu+1,b}f(z)|< M.
\end{equation*}
\begin{proof}
Let $\phi(\alpha,\beta,\gamma,\delta;z)=\beta-\alpha$ and $\Omega=h(\mathbb{U}),$
where
$h(z)=\frac{Mz}{|b+1|}\quad (M>0)$.
Use Corollary \ref{c3}, we need to show that $\phi \in \Phi_{J}[\Omega,M],$ that is, the admissibility condition (\ref{p16}) is satisfied. This follows since
\begin{equation*}
|\phi(v,w,x,y;z)|=\left|\frac{(k-1)Me^{i\theta}}{b+1}\right|\geq \frac{M}{|b+1|},
\end{equation*}
whenever $z \in \mathbb{U},~~ \theta \in \mathbb{R}$ and $k \geq 2.$ The required result follows from Corollary \ref{c3}.
\end{proof}
\end{cor}
\begin{definition}\label{d3}
{\rm Let $\Omega$ be a set in $\mathbb{C}, q \in \mathbb{Q}_{1} \cap \mathcal{H}_{1}.$ The class of admissible functions $\Phi_{J,1}[\Omega,q]$ consists of those functions $\phi:\mathbb{C}^{4}\times \mathbb{U}\longrightarrow \mathbb{C}$ that satisfy the following admissibility condition
\begin{equation*}
 \phi(\alpha,\beta,\gamma,\delta;z) \not\in \Omega
 \end{equation*}
whenever
\begin{equation*}
\alpha=q(\zeta),\beta=\frac{k\zeta q^{\prime}(\zeta)+(b+1)q(\zeta)}{b+1}, 
\end{equation*}
\begin{equation*}
\Re\left(\frac{(b+1)(\gamma-\alpha)}{\beta-\alpha}-2(1+b)\right)\geq k \Re\left(\frac{\zeta q^{\prime\prime}(\zeta)}{q^{\prime}(\zeta)}+1\right),
\end{equation*}
and 
\begin{eqnarray*}
&&\Re\left(\frac{\delta(1+b)^{2}-3\gamma(b+2)(b+1)+3\alpha(b+2)(b+1)-(1+b)^{2}\alpha}{\beta-\alpha}+3b^{2}+12b+11\right)\geq k^{2} \Re\left(\frac{\zeta^{2} q^{\prime\prime\prime}(\zeta)}{q^{\prime}(\zeta)}\right),
\end{eqnarray*}
where $z \in \mathbb{U},\zeta \in \partial \mathbb{U} \setminus E(q),$ and $k\geq 2.$}
\end{definition}
\begin{thm}\label{t9}
Let $\phi \in \Phi_{J,1}[\Omega,q].$ If the function $f \in \mathcal{A}$ and $q \in \mathbb{Q}_{1}$ satisfy the following conditions:
\begin{equation}\label{p17}
\Re\left(\frac{\zeta q^{\prime\prime}(\zeta)}{q^{\prime}(\zeta)}\right)\geq 0\quad \left|\frac{J_{\mu,b}f(z)}{zq^{\prime}(\zeta)}\right|\leq k,
\end{equation}
and 
\begin{equation}\label{p18}
\left\{\phi \left(\frac{J_{\mu+1,b}f(z)}{z},\frac{J_{\mu,b}f(z)}{z},\frac{J_{\mu-1,b}f(z)}{z},\frac{J_{\mu-2,b}f(z)}{z};z \right):z \in \mathbb{U}\right\}\subset \Omega,
\end{equation}
then 
\begin{equation*}
\frac{J_{\mu+1,b}f(z)}{z} \prec q(z)\qquad(z \in \mathbb{U}).
\end{equation*}

\begin{proof}
Define the analytic function $p(z)$ in $\mathbb{U}$ by
\begin{equation}\label{p19}
p(z)=\frac{J_{\mu+1,b}f(z)}{z}.
\end{equation}
From equation (\ref{g5}) and (\ref{p19}), we have
\begin{equation}\label{p20}
\frac{J_{\mu,b}f(z)}{z}=\frac{zp^{\prime}(z)+(b+1)p(z)}{b+1}.
\end{equation}
By similar argument, yields 
\begin{equation}\label{p21}
\frac{J_{\mu-1,b}f(z)}{z}=\frac{z^{2}p^{\prime\prime}(z)+zp^{\prime}(z)(3+2b)+p(z)(1+b)^{2}}{(b+1)^{2}}
\end{equation}
and
\begin{equation}\label{p22}
\frac{J_{\mu-2,b}f(z)}{z}=\frac{z^{3}p^{\prime\prime\prime}(z)+3(b+2)z^{2}p^{\prime\prime}(z)+(3b^{2}+9b+7)zp^{\prime}(z)+p(z)(b+1)^{3}}{(b+1)^{3}}.
\end{equation}
Define the transformation from $\mathbb{C}^{4}$ to $\mathbb{C}$ by
\begin{equation*}
\alpha(r,s,t,u)=r,\qquad \beta(r,s,t,u)=\frac{s+(b+1)r}{(b+1)},
\end{equation*}
\begin{equation}\label{p23}
\gamma(r,s,t,u)=\frac{t+(3+2b)s+(b+1)^{2}r}{(b+1)^{2}},
\end{equation}
and
\begin{equation}\label{p24}
\delta(r,s,t,u)=\frac{u+3(b+2)t+(3b^{2}+9b+7)s+(b+1)^{3}r}{(b+1)^{3}}.
\end{equation}
Let
\begin{multline}\label{p25}
\psi(r,s,t,u)=\phi(\alpha,\beta,\gamma,\delta;z)=\phi \bigg(r,\frac{s+(1+b)r}{(1+b)},\frac{t+(3+2b)s+(b+1)^{2}r}{(b+1)^{2}},\\ \frac{u+3(b+2)t+(3b^{2}+9b+7)s+(b+1)^{3}r}{(b+1)^{3}};z\bigg)
\end{multline}
The proof will make use of  Lemma \ref{t4}. Using equations (\ref{p19}) to (\ref{p22}),  and from (\ref{p25}), we have
\begin{equation}\label{p26}
\psi \left(p(z),zp^{\prime}(z),z^{2}p^{\prime\prime}(z),z^{3}p^{\prime\prime\prime}(z);z\right)=\phi \left(\frac{J_{\mu+1,b}f(z)}{z},\frac{J_{\mu,b}f(z)}{z},\frac{J_{\mu-1,b}f(z)}{z},\frac{J_{\mu-2,b}f(z)}{z};z \right).
\end{equation} 
Hence,(\ref{p18}) becomes
\begin{equation*}
\psi \left(p(z),zp^{\prime}(z),z^{2}p^{\prime\prime}(z),z^{3}p^{\prime\prime\prime}(z);z\right) \in \Omega.
\end{equation*}
Note that
\begin{equation*}
\frac{t}{s}+1=\frac{(b+1)(\gamma-\alpha)}{\beta-\alpha}-2(1+b)
\end{equation*}
and
\begin{eqnarray*}
\frac{u}{s}=\frac{\delta(1+b)^{2}-3\gamma(b+2)(b+1)+3\alpha(b+2)(b+1)-(1+b)^{2}\alpha}{\beta-\alpha}+3b^{2}+12b+11.
\end{eqnarray*}
Thus, the admissibility condition for $\phi \in \Phi_{J,1}[\Omega,q]$ in Definition \ref{d3} is equivalent to the admissibility condition for $\psi \in \Psi_{2}[\Omega,q]$ as given in Definition \ref{d2a} with $n=2.$ Therefore, by using (\ref{p17}) and Lemma \ref{t4}, we have
\begin{equation*}
\frac{J_{\mu+1,b}f(z)}{z}\prec q(z).
\end{equation*}
This completes the proof of the theorem.
\end{proof}
\end{thm}
If $\Omega \neq \mathbb{C}$ is a simply connected domain, then $\Omega=h(\mathbb{U})$ for some conformal mapping $h(z)$ of $\mathbb{U}$ onto $\Omega.$ In this case, the class $\Phi_{J,1}[h(\mathbb{U}),q]$ is written  as $\Phi_{J,1}[h,q].$ This follows immediate  consequence of Theorem  \ref{t9}.
\begin{thm}\label{t10}
Let $\phi \in \Phi_{J,1}[h,q].$ If the function $f \in \mathcal{A}$ and $q \in \mathbb{Q}_{1}$ satisfy the following conditions:
\begin{equation}\label{p27}
\Re\left(\frac{\zeta q^{\prime\prime}(\zeta)}{q^{\prime}(\zeta)}\right)\geq 0,\quad \left|\frac{J_{\mu,b}f(z)}{zq^{\prime}(\zeta)}\right|\leq k,
\end{equation}
and 
\begin{equation}\label{p28}
\phi \left(\frac{J_{\mu+1,b}f(z)}{z},\frac{J_{\mu,b}f(z)}{z},\frac{J_{\mu-1,b}f(z)}{z},\frac{J_{\mu-2,b}f(z)}{z};z \right)\prec h(z),
\end{equation}
then 
\begin{equation*}
\frac{J_{\mu+1,b}f(z)}{z} \prec q(z)\qquad(z \in \mathbb{U}).
\end{equation*}
\end{thm}
In view of Definition \ref{d3} and in the special case $q(z)=Mz,\quad M>0,$ the class of admissible functions $\Phi_{J,1}[\Omega,q],$ denoted by  $\Phi_{J,1}[\Omega,M],$ is expressed as follows.
\begin{definition}\label{d40}
{\rm Let $\Omega$ be a set in $\mathbb{C}$ and $M>0.$ The class of admissible function $\Phi_{J,1}[\Omega,M]$ consists of those functions $\phi:\mathbb{C}^{4}\times \mathbb{U}\longrightarrow\mathbb{C}$ such that
\begin{multline}\label{p29}
\phi \bigg(Me^{i\theta},\frac{(k+1+b)Me^{i\theta}}{1+b},\frac{L+[(3+2b)k+(b+1)^{2}]Me^{i\theta}}{(b+1)^{2}},\\ \frac{N+3(b+2)L+[(3b^{2}+9b+7)k+(b+1)^{3}]Me^{i\theta}}{(b+1)^{3}};z\bigg)\notin \Omega.
\end{multline}
whenever $z \in \mathbb{U},\Re (Le^{-i\theta})\geq (k-1)kM,$ and $\Re (Ne^{-i\theta})\geq 0$ for all $\theta \in \mathbb{R}$ and $k \geq 2.$ }
\end{definition}
\begin{cor}\label{c5}
Let $\phi \in \Phi_{J,1}[\Omega,M].$ If the function $f \in \mathcal{A}$ satisfies 
\begin{equation*}
\left|\frac{J_{\mu,b}f(z)}{z}\right|\leq kM\qquad(z \in \mathbb{U},k\geq 2 ; M>0),
\end{equation*}
and 
\begin{equation*}
\phi \left(\frac{J_{\mu+1,b}f(z)}{z},\frac{J_{\mu,b}f(z)}{z},\frac{J_{\mu-1,b}f(z)}{z},\frac{J_{\mu-2,b}f(z)}{z};z \right)\in \Omega,
\end{equation*}
then 
\begin{equation*}
\left|\frac{J_{\mu+1,b}f(z)}{z}\right|<M.
\end{equation*}
\end{cor}
 In this special case $\Omega=q(\mathbb{U})=\{w:|w|<M\},$ the class $\Phi_{J,1}[\Omega,M]$ is simply denoted by $\Phi_{J,1}[M].$ Corollary \ref{c5} can now be written in the following form.
\begin{cor}\label{c6}
 Let $\phi \in \Phi_{J,1}[M].$ If the function $f \in \mathcal{A}$ satisfies 
\begin{equation*}
\left|\frac{J_{\mu,b}f(z)}{z}\right|\leq kM\qquad(z \in \mathbb{U},k\geq 2; M>0),
\end{equation*}
and 
\begin{equation*}
\left|\phi\left(\frac{J_{\mu+1,b}f(z)}{z},\frac{J_{\mu,b}f(z)}{z},\frac{J_{\mu-1,b}f(z)}{z},\frac{J_{\mu-2,b}f(z)}{z}\right);z \right|<M,
\end{equation*}
then
\begin{equation*}
\left|\frac{J_{\mu+1,b}f(z)}{z}\right|< M.
\end{equation*}
\end{cor}
\begin{definition}\label{d50}
{\rm Let $\Omega$  be a set in $\mathbb{C}$ and $q \in \mathbb{Q}_{1} \cap \mathcal{H}_{1}.$ The class of admissible functions $\Phi_{J,2}[\Omega,q]$ consists of those functions $\phi: \mathbb{C}^{4} \times \mathbb{U}\longrightarrow \mathbb{C}$ that satisfy the following admissibility condition
\begin{equation*}
\phi(\alpha,\beta,\gamma,\delta;z) \not\in \Omega
\end{equation*}
whenever
\begin{equation*}
\alpha=q(\zeta), \beta=\frac{1}{(b+1)}\left(\frac{k\zeta q^{\prime}(\zeta)}{q(\zeta)}+(b+1) q(\zeta)\right),
\end{equation*}
\begin{equation*}
\Re\left(\frac{(1+b)(\beta \gamma+2\alpha^{2}-3\alpha \beta)}{(\beta-\alpha)}\right)\geq k \Re\left(\frac{\zeta q^{\prime\prime}(\zeta)}{q^{\prime}(\zeta)}+1\right),
\end{equation*}
and
\begin{multline*}
\Re \bigg[(\delta-\gamma)(1+b)^{2}\beta \gamma-(1+b)^{2}(\gamma-\beta)\beta(1-\beta-\gamma+3\alpha)-3(b+1)(\gamma-\beta)\beta+2(\beta-\alpha)+3(1+b)\alpha(\beta-\alpha)\\ \qquad+(\beta-\alpha)^{2}(1+b)((\beta-\alpha)(1+b)-3-4(1+b)\alpha)+\alpha^{2}(1+b)^{2}(\beta-\alpha)
\bigg]\times (\beta-\alpha)^{-1}\geq k^{2}\Re \bigg(\frac{{\zeta}^{2}q^{\prime\prime\prime}(\zeta)}{q^{\prime}(\zeta)}\bigg),
\end{multline*}
where $z \in \mathbb{U},\zeta \in \partial \mathbb{U} \setminus E(q)$ and $k\geq 2.$}
\end{definition}
\begin{thm}\label{t11}
Let $\phi \in \Phi_{J,2}[\Omega,q].$ If the function $f \in \mathcal{A}$ and $q \in \mathbb{Q}_{1}$ satisfy the following conditions
\begin{equation}\label{e30}
\Re\left(\frac{\zeta q^{\prime\prime}(\zeta)}{q^{\prime}(\zeta)}\right)\geq 0, \quad\left| \frac{J_{\mu-1,b}f(z)}{J_{\mu,b}f(z)q^{\prime}(\zeta)}\right|\leq k,
\end{equation} 
\begin{equation}\label{e31}
\left\{\phi \left(\frac{J_{\mu,b}f(z)}{J_{\mu+1,b}f(z)},\frac{J_{\mu-1,b}f(z)}{J_{\mu,b}f(z)},\frac{J_{\mu-2,b}f(z)}{J_{\mu-1,b}f(z)},\frac{J_{\mu-3,b}f(z)}{J_{\mu-2,b}f(z)};z\right):z \in \mathbb{U}\right\}\subset \Omega,
\end{equation}
then
\begin{equation*}
\frac{J_{\mu,b}f(z)}{J_{\mu+1,b}f(z)}\prec q(z)\qquad(z \in \mathbb{U}).
\end{equation*}
\begin{proof}
Define the analytic function $p(z)$ in $\mathbb{U}$ by
\begin{equation}\label{e32}
p(z)=\frac{J_{\mu,b}f(z)}{J_{\mu+1,b}f(z)}.
\end{equation}
From equation (\ref{g5}) and (\ref{e32}), we have
\begin{equation}\label{e33}
\frac{J_{\mu-1,b}f(z)}{J_{\mu,b}f(z)}=\frac{1}{(b+1)}\left[\frac{zp^{\prime}(z)}{p(z)}+(b+1)p(z) \right]:=\frac{A}{b+1}.
\end{equation}
By similar argument yields,
\begin{equation}\label{e34}
\frac{J_{\mu-2,b}f(z)}{J_{\mu-1,b}f(z)}:=\frac{B}{b+1}
\end{equation}
and
\begin{eqnarray}\label{e35}
\frac{J_{\mu-3,b}f(z)}{J_{\mu-2,b}f(z)}=\frac{1}{b+1}\left[B+B^{-1}(C+A^{-1}D-A^{-2}C^{2})\right].
\end{eqnarray}
Where
\begin{align*}
B&:=(b+1)p(z)+\frac{zp^{\prime}(z)}{p(z)}+\dfrac{\frac{z^{2}p^{\prime\prime}(z)}{p(z)}+\frac{zp^{\prime}(z)}{p(z)}-\left(\frac{zp^{\prime}(z)}{p(z)}\right)^{2}+(b+1)zp^{\prime}(z)}{\frac{zp^{\prime}(z)}{p(z)}+(b+1)p(z)}\\ 
C&:=\frac{z^{2}p^{\prime\prime}(z)}{p(z)}+\frac{zp^{\prime}(z)}{p(z)}-\left(\frac{zp^{\prime}(z)}{p(z)}\right)^{2}+(b+1)zp^{\prime}(z)\\
D&:=\frac{3z^{2}p^{\prime\prime}(z)}{p(z)}+  \frac{z^{3}p^{\prime\prime\prime}(z)}{p(z)}+\frac{zp^{\prime}(z)}{p(z)}-3\left(\frac{zp^{\prime}(z)}{p(z)}\right)^{2}-\frac{3z^{3}p^{\prime}(z)p^{\prime\prime}(z)}{p^{2}(z)}+2\left(\frac{zp^{\prime}(z)}{p(z)}\right)^{3} \\\qquad \qquad &+(b+1)zp^{\prime}(z)+(b+1)z^{2}p^{\prime\prime}(z).
\end{align*}
Define the transformation from $\mathbb{C}^{4}$ to $\mathbb{C}$ by
\begin{eqnarray*}
\alpha(r,s,t,u)=r,\qquad \beta(r,s,t,u)=\frac{1}{b+1}\left[\dfrac{s}{r}+(b+1)r \right]:=\frac{E}{b+1},
\end{eqnarray*}
\begin{eqnarray}\label{e36}
\gamma(r,s,t,u)=\frac{1}{b+1}\left[\frac{s}{r}+(b+1)r+\frac{\frac{t}{r}+\frac{s}{r}-(\frac{s}{r})^{2}+(b+1)s}{\frac{s}{r}+(b+1)r} \right]:=\frac{F}{b+1},
\end{eqnarray}
and
\begin{eqnarray}\label{e37}
\delta(r,s,t,u)=\frac{1}{b+1}\left[F+F^{-1}(L+E^{-1}H-E^{-2}L^{2})\right].
\end{eqnarray}
Where
\begin{align*}
L&:=(1+b)s+\frac{t}{r}+\frac{s}{r}-\left(\frac{s}{r}\right)^{2}\\
H&:=\frac{3t}{r}+\frac{u}{r}+\frac{s}{r}-3\left(\frac{s}{r}\right)^{2}-3\left(\frac{st}{r^{2}}\right)+2\left(\frac{s}{r}\right)^{3}+(1+b)s+(1+b)t.
\end{align*}
Let
\begin{multline}\label{e38}
\psi(r,s,t,u)=\phi(\alpha,\beta,\gamma,\delta;z)=\phi \left(r,\frac{E}{b+1},\frac{F}{b+1},\frac{F+F^{-1}(L+E^{-1}H-E^{-2}L^{2})}{b+1}\right).
\end{multline}

The proof will make use of  Lemma  \ref{t4}. Using equations (\ref{e32}) to (\ref{e35}), and from (\ref{e38}), we have
\begin{equation}\label{e39}
\psi \left(p(z),zp^{\prime}(z),z^{2}p^{\prime\prime}(z),z^{3}p^{\prime\prime\prime}(z);z\right)=\phi \left(\frac{J_{\mu,b}f(z)}{J_{\mu+1,b}f(z)},\frac{J_{\mu-1,b}f(z)}{J_{\mu,b}f(z)},\frac{J_{\mu-2,b}f(z)}{J_{\mu-1,b}f(z)},\frac{J_{\mu-3,b}f(z)}{J_{\mu-2,b}f(z)};z\right).
\end{equation} 
Hence,(\ref{e31}) becomes
\begin{equation*}
\psi \left(p(z),zp^{\prime}(z),z^{2}p^{\prime\prime}(z),z^{3}p^{\prime\prime\prime}(z);z\right) \in \Omega.
\end{equation*}
Note that
\begin{equation*}
\frac{t}{s}+1=\left(\frac{(1+b)(\beta \gamma+2\alpha^{2}-3\alpha \beta)}{(\beta-\alpha)}\right)
\end{equation*}
and
\begin{multline*}
\frac{u}{s}= \bigg[(\delta-\gamma)(1+b)^{2}\beta \gamma-(1+b)^{2}(\gamma-\beta)\beta(1-\beta-\gamma+3\alpha)-3(b+1)(\gamma-\beta)\beta+2(\beta-\alpha)+3(1+b)\alpha(\beta-\alpha)\\+(\beta-\alpha)^{2}(1+b)((\beta-\alpha)(1+b)-3-4(1+b)\alpha)+\alpha^{2}(1+b)^{2}(\beta-\alpha)
\bigg]\times (\beta-\alpha)^{-1}.
\end{multline*}
Thus, the admissibility condition for $\phi \in \Phi_{J,2}[\Omega,q]$ in Definition \ref{d50} is equivalent to the admissibility condition for $\psi \in \Psi_{2}[\Omega,q]$ as given in Definition \ref{d2a} with $n=2.$ Therefore, by using (\ref{e30}) and Lemma \ref{t4}, we have
\begin{equation*}
\frac{J_{\mu,b}f(z)}{J_{\mu+1,b}f(z)}\prec q(z).
\end{equation*}
\end{proof}
\end{thm}
If $\Omega \neq \mathbb{C}$ is a simply connected domain, then $\Omega=h(\mathbb{U})$ for some conformal mapping $h(z)$ of $\mathbb{U}$ onto $\Omega.$ In this case, the class $\Phi_{J,1}[h(\mathbb{U}),q]$ is written  as $\Phi_{J,2}[h,q]$. This follows immediate consequence of Theorem  \ref{t11}.
\begin{thm}\label{t12}
Let $\phi \in \Phi_{J,2}[h,q].$ If the function $f \in \mathcal{A}$ and $q \in \mathbb{Q}_{1}$ satisfy the following conditions (\ref{e30}) and 
\begin{equation}\label{e40}
\phi \left(\frac{J_{\mu,b}f(z)}{J_{\mu+1,b}f(z)},\frac{J_{\mu-1,b}f(z)}{J_{\mu,b}f(z)},\frac{J_{\mu-2,b}f(z)}{J_{\mu-1,b}f(z)},\frac{J_{\mu-3,b}f(z)}{J_{\mu-2,b}f(z)};z\right)\prec h(z),
\end{equation}
then 
\begin{equation*}
\frac{J_{\mu,b}f(z)}{J_{\mu+1,b}f(z)} \prec q(z)\qquad(z \in \mathbb{U}).
\end{equation*}
\end{thm}
\section{\bf Results Related to Third Order Superordination}

In this section, the third-order differential superordination theorems for the operator $J_{\mu,b}f(z)$ defined in (\ref{g4}) is investigated. For the purpose, we considered the following admissible functions.
\begin{definition}\label{d7}
{\rm Let $\Omega$ be a set in $\mathbb{C}$ and $q \in \mathcal{H}_{0}$ with $q^{\prime}(z) \neq 0.$ The class of admissible function $\Phi_{J}^{\prime}[\Omega,q]$ consists of those functions $\phi:\mathbb{C}^{4}\times \mathcal{\overline{U}}\longrightarrow \mathbb{C}$ that satisfy the following admissibility conditions:
\begin{equation*}
\phi(\alpha,\beta,\gamma,\delta;\zeta) \in \Omega
\end{equation*}
whenever
\begin{equation*}
\alpha=q(z), \beta=\frac{z q^{\prime}(z)+m b q(z)}{m(b+1)},
\end{equation*}
\begin{equation*}
\Re\left(\frac{\gamma(b+1)^{2}-b^{2}\alpha}{(\beta(b+1)-b\alpha)}-2b\right)\leq \frac{1}{m} \Re\left(\frac{z q^{\prime\prime}(z)}{q^{\prime}(z)}+1\right),
\end{equation*}
and 
\begin{equation*}
\Re\left(\frac{\delta(b+1)^{3}-\gamma(b+1)^{2}(3b+3)+b^{2}\alpha(3+2b)}{(b(\beta-\alpha)+\beta)}+3b^{2}+6b+2\right)\leq \frac{1}{m^{2}} \Re \left(\frac{z^{2} q^{\prime\prime\prime}(z)}{q^{\prime}(z)}\right),
\end{equation*}
where $z \in \mathbb{U},\zeta \in \partial \mathbb{U},$ and $m\geq 2.$} \end{definition}
\begin{thm}\label{t13}
Let $\phi \in \Phi_{J}^{\prime}[\Omega,q].$ If the function $f \in \mathcal{A}$ and $J_{\mu+1,b}f(z)\in \mathbb{Q}_{0}$ and $q \in \mathcal{H}_{0}$ with $q^{\prime}(z)\neq 0$ satisfy the following conditions:
\begin{equation}\label{p30}
\Re\left(\frac{z q^{\prime\prime}(z)}{q^{\prime}(z)}\right)\geq 0,\quad \left|\frac{J_{\mu,b}f(z)}{q^{\prime}(z)}\right|\leq m,
\end{equation}
and 
\begin{equation*}
\phi(J_{\mu+1,b}f(z),J_{\mu,b}f(z),J_{\mu-1,b}f(z),J_{\mu-2,b}f(z);z)
\end{equation*}
is univalent in $\mathbb{U}$,then 
 \begin{equation}\label{p31}
\Omega \subset \{\phi(J_{\mu+1,b}f(z),J_{\mu,b}f(z),J_{\mu-1,b}f(z),J_{\mu-2,b}f(z)
;z):z \in \mathbb{U}\},
 \end{equation}
 implies that
 \begin{equation*}
q(z) \prec J_{\mu+1,b}f(z)\qquad(z \in \mathbb{U}).
\end{equation*}
\begin{proof}
Let the function $p(z)$ be defined by (\ref{p5}) and $\psi$ by  (\ref{p11}). Since $\phi \in \Phi_{J}^{\prime}[\Omega,q].$ From (\ref{p12}) and (\ref{p31}) yield
\begin{equation*}
\Omega \subset \{\psi \left(p(z),zp^{\prime}(z),z^{2}p^{\prime\prime}(z),z^{3}p^{\prime\prime\prime}(z);z\right):z \in \mathbb{U}\}.
\end{equation*}
From (\ref{p9}) and (\ref{p10}), we see that  the admissibility condition for $\phi \in \Phi_{J}^{\prime}[\Omega,q]$ in Definition \ref {d7} is equivalent to the admissibility condition for $\psi \in \Psi_{2}[\Omega,q]$ as given in Definition \ref {d6} with $n=2.$ Hence $\psi \in \Psi_{2}^{\prime}[\Omega,q]$  and by using (\ref{p31}) and Lemma \ref{t5}, we have
\begin{equation*}
q(z) \prec J_{\mu+1,b}f(z).
\end{equation*}
\end{proof}
\end{thm}
If $\Omega \neq \mathbb{C}$ is a simply connected domain, then $\Omega=h(\mathbb{U})$ for some conformal mapping $h(z)$ of $\mathbb{U}$ onto $\Omega.$ In this case, the class $\Phi_{J}^{\prime}[h(\mathbb{U}),q]$ is written  as $\Phi_{J}^{\prime}[h,q].$ This follows an immediate  consequence of Theorem  \ref{t13}.

\begin{thm}\label{t14}
Let $\phi \in \Phi_{J}^{\prime}[h,q]$ and $h$ be analytic in $\mathbb{U}$. If the function $f \in \mathcal{A}$ and $J_{\mu+1,b}f(z) \in \mathbb{Q}_{0}$ and $q \in  \mathcal{H}_{0}$  with $q ^{\prime}(z)\neq 0$ satisfy the following conditions (\ref{p30})
and 
\begin{equation*}
\phi(J_{\mu+1,b}f(z),J_{\mu,b}f(z),J_{\mu-1,b}f(z),J_{\mu-2,b}f(z)
;z),
\end{equation*}
is univalent in $\mathbb{U},$ then
\begin{equation}\label{p32}
h(z) \prec\phi(J_{\mu+1,b}f(z),J_{\mu,b}f(z),J_{\mu-1,b}f(z),J_{\mu-2,b}f(z);z),
\end{equation}
implies that 
\begin{equation*}
q(z) \prec J_{\mu+1,b}f(z)\qquad(z \in \mathbb{U}).
\end{equation*}
\end{thm}
Theorems \ref{t13} and \ref{t14} can only be used to obtain subordination of the third-order differential superordination of the forms (\ref{p31}) or (\ref{p32}). The following theorem proves the existence of the best subordinant of (\ref{p32}) for a suitable  $\phi.$
\begin{thm}\label{t15}
Let the function $h$ be univalent in $\mathbb{U}$ and let $\phi : \mathbb{C}^{4}\times \mathcal{\overline{U}}\longrightarrow \mathbb{C}$ and $\psi$ be given by (\ref{p11}). Suppose that the differential equation
\begin{equation}\label{p33}
\psi(q(z),zq^{\prime}(z),z^{2}q^{\prime\prime}(z),z^{3}
q^{\prime\prime\prime}(z);z)=h(z)
\end{equation}
has a solution $q(z) \in \mathbb{Q}_{0}.$ If the functions $f \in \mathcal{A}, J_{\mu+1,b}f(z) \in \mathbb{Q}_{0}$ and $q \in \mathcal{H}_{0}$ with $q^{\prime}(z) \neq 0,$ which satisfy  the following condition (\ref{p30}) and 
\begin{equation*}
\phi(J_{\mu+1,b}f(z),J_{\mu,b}f(z),J_{\mu-1,b}f(z),J_{\mu-2,b}f(z)
;z)
\end{equation*}
is analytic in $\mathbb{U},$ then
\begin{equation*}
h(z) \prec \phi(J_{\mu+1,b}f(z),J_{\mu,b}f(z),J_{\mu-1,b}f(z),J_{\mu-2,b}f(z)
;z)
\end{equation*}
implies that
\begin{equation*}
q(z) \prec J_{\mu+1,b}f(z)\qquad(z \in \mathbb{U})
\end{equation*}
and $q(z)$ is the best dominant.
\begin{proof}
In view of Theorem \ref{t13} and Theorem \ref{t14}, we deduce that $q$ is a subordinant of (\ref{p32}). Since $q$ satisfies (\ref{p33}), it is also a solution of (\ref{p32}) and therefore $q$ will be subordinated by all subordinants. Hence $q$ is the best subordinant. 
\end{proof}
\end{thm}
\begin{definition}\label{d8}
{\rm Let $\Omega$ be a set in $\mathbb{C}, q \in \mathcal{H}_{1}$ with $q^{\prime}(z) \neq 0$.  The class of admissible functions $\Phi_{J,1}^{\prime}[\Omega,q]$ consists of those functions $\phi:\mathbb{C}^{4}\times \mathcal{\overline{U}}\longrightarrow \mathbb{C}$ that satisfy the following admissibility condition
\begin{equation*}
 \phi(\alpha,\beta,\gamma,\delta;\zeta) \in \Omega
\end{equation*}
whenever
\begin{equation*}
\alpha=q(z),\beta=\frac{z q^{\prime}(z)+(1+b)mq(z)}{(1+b)m}, 
\end{equation*}
\begin{equation*}
\Re\left(\frac{(b+1)(\gamma-\alpha)}{\beta-\alpha}-2(1+b)\right)\leq \frac{1}{m} \Re\left(\frac{z q^{\prime\prime}(z)}{q^{\prime}(z)}+1\right)
\end{equation*}
and 
\begin{eqnarray*}
&&\Re\left(\frac{\delta(1+b)^{2}-3\gamma(b+2)(b+1)+3\alpha(b+2)(b+1)-(1+b)^{2}\alpha}{\beta-\alpha}+3b^{2}+12b+11\right)\leq \frac{1}{m^{2}} \Re\left(\frac{z^{2} q^{\prime\prime\prime}(z)}{q^{\prime}(z)}\right),
\end{eqnarray*}
where $z \in \mathbb{U},\zeta \in \partial \mathbb{U},$ and $m\geq 2.$}
\end{definition}
\begin{thm}\label{t16}
Let $\phi \in \Phi_{J,1}^{\prime}[\Omega,q].$ If the function $f \in \mathcal{A}, \frac{J_{\mu,b}f(z)}{z} \in \mathbb{Q}_{1}$ and $q \in \mathcal{H}_{1}$ with $q^{\prime}(z) \neq 0$ satisfy the following conditions:
\begin{equation}\label{p34}
\Re\left(\frac{z q^{\prime\prime}(z)}{q^{\prime}(z)}\right)\geq 0,\left|\frac{J_{\mu,b}f(z)}{zq^{\prime}(z)}\right|\leq m,
\end{equation}
and 
\begin{equation*}
\phi \left(\frac{J_{\mu+1,b}f(z)}{z},\frac{J_{\mu,b}f(z)}{z},\frac{J_{\mu-1,b}f(z)}{z},\frac{J_{\mu-2,b}f(z)}{z};z \right),
\end{equation*}
is univalent  in $\mathbb{U},$
then 
\begin{equation}\label{p35}
\Omega \subset \left\{\phi \left(\frac{J_{\mu+1,b}f(z)}{z},\frac{J_{\mu,b}f(z)}{z},\frac{J_{\mu-1,b}f(z)}{z},\frac{J_{\mu-2,b}f(z)}{z};z \right): z \in \mathbb{U} \right\},
\end{equation}
implies that
\begin{equation*}
q(z) \prec \frac{J_{\mu+1,b}f(z)}{z}\qquad(z \in \mathbb{U}).
\end{equation*}
\begin{proof}
Let the function $p(z)$ be defined by (\ref{p19}) and $\psi$ by (\ref{p25}). Since $\phi \in \Phi _{J,1}^{\prime}[\Omega,q],$ (\ref{p26}) and (\ref{p35}) yield
\begin{equation*}
\Omega \subset \left\{\psi \left(p(z),zp^{\prime}(z),z^{2}p^{\prime\prime}(z),z^{3}p^{\prime\prime\prime}(z);z\right):z \in \mathbb{U}\right\}.
\end{equation*}
From equations (\ref{p23}) and (\ref{p24}), we see that the admissible condition for $\phi \in \Phi_{J,1}^{\prime}[\Omega,q]$ in Definition \ref{d8} is equivalent to the admissible condition for $\psi$ as given in Definition \ref{d6} with $n=2.$ Hence $\psi \in \Psi _{2}^{\prime}[\Omega,q]$ and by using (\ref{p34}) and Lemma \ref{t5}, we have
\begin{equation*}
q(z)\prec \frac{J_{\mu+1,b}f(z)}{z}.
\end{equation*} 
\end{proof}
\end{thm}
If $\Omega \neq \mathbb{C}$ is a simply connected domain, then $\Omega=h(\mathbb{U})$ for some conformal mapping $h(z)$ of $\mathbb{U}$ onto $\Omega.$ In this case, the class $\Phi_{J,1}^{\prime}[h(\mathbb{U}),q]$ is written  as $\Phi_{J,1}^{\prime}[h,q].$ This follows an immediate  consequence of Theorem  \ref{t16}.
\begin{thm}\label{t17}
Let $\phi \in \Phi_{J,1}^{\prime}[h,q]$ and $h$ be analytic in $\mathbb{U}$. If the function $f \in \mathcal{A}$ and $q \in \mathcal{H}_{1}$ with $q^{\prime}(z)\neq 0$ satisfy the following conditions (\ref{p34})
and 
\begin{equation*}
\phi \left(\frac{J_{\mu+1,b}f(z)}{z},\frac{J_{\mu,b}f(z)}{z},\frac{J_{\mu-1,b}f(z)}{z},\frac{J_{\mu-2,b}f(z)}{z};z \right),
\end{equation*}
is univalent in $\mathbb{U},$
then
\begin{equation*}
h(z)\prec\phi \left(\frac{J_{\mu+1,b}f(z)}{z},\frac{J_{\mu,b}f(z)}{z},\frac{J_{\mu-1,b}f(z)}{z},\frac{J_{\mu-2,b}f(z)}{z};z \right),
\end{equation*}
implies that
\begin{equation*}
q(z) \prec\frac{J_{\mu+1,b}f(z)}{z}\qquad(z \in \mathbb{U}).
\end{equation*}
\end{thm}

\begin{definition}\label{d9}
{\rm Let $\Omega$  be a set in $\mathbb{C}$ and $q \in \mathcal{H}_{1}$ with $q^{\prime}(z)\neq 0.$ The class of admissible functions $\Phi_{J,2}^{\prime}[\Omega,q]$ consists of those functions $\phi: \mathbb{C}^{4} \times \mathcal{\overline{U}}\longrightarrow \mathbb{C}$ that satisfy the following admissibility conditions
\begin{equation*}
\phi(\alpha,\beta,\gamma,\delta;\zeta) \in \Omega
\end{equation*}
whenever
\begin{equation*}
\alpha=q(z),\beta=\frac{1}{b+1}\left(\frac{z q^{\prime}(z)}{mq(z)}+(b+1)q(z)\right),
\end{equation*}
\begin{equation*}
\Re \left(\frac{(1+b)(\beta \gamma+2\alpha^{2}-3\alpha \beta)}{(\beta-\alpha)}\right)\leq \frac{1}{m}\Re \left(\frac{z q^{\prime\prime}(z)}{q^{\prime}(z)}+1\right),
\end{equation*}and
\begin{multline*}
\Re \bigg[(\delta-\gamma)(1+b)^{2}\beta \gamma-(1+b)^{2}(\gamma-\beta)\beta(1-\beta-\gamma+3\alpha)-3(b+1)(\gamma-\beta)\beta+2(\beta-\alpha)+3(1+b)\alpha(\beta-\alpha)\\ \qquad+(\beta-\alpha)^{2}(1+b)((\beta-\alpha)(1+b)-3-4(1+b)\alpha)+\alpha^{2}(1+b)^{2}(\beta-\alpha)
\bigg]\times (\beta-\alpha)^{-1} \leq \frac{1}{m^{2}}\Re\left(\frac{{z}^{2}q^{\prime\prime\prime}(z)}{q^{\prime}(z)}\right),
\end{multline*}
where $z \in \mathbb{U},\zeta \in \partial \mathbb{U}$ and $m\geq 2.$}
\end{definition}
\begin{thm}\label{t18}
Let $\phi \in \Phi_{J,2}^{\prime}[\Omega,q].$ If the function $f \in \mathcal{A}$ and $\frac{J_{\mu,b}f(z)}{J_{\mu+1,b}f(z)} \in \mathbb{Q}_{1}$ and $q \in \mathcal{H}_{1}$ with $q^{\prime}(z)\neq 0$ satisfy the following conditions
\begin{equation}\label{p36}
\Re\left(\frac{z q^{\prime\prime}(z)}{q^{\prime}(z)}\right)\geq 0,\left| \frac{J_{\mu-1,b}f(z)}{J_{\mu,b}f(z)q^{\prime}(z)}\right|\leq m,
\end{equation} 
and
\begin{equation*}
\phi \left(\frac{J_{\mu,b}f(z)}{J_{\mu+1,b}f(z)},\frac{J_{\mu-1,b}f(z)}{J_{\mu,b}f(z)},\frac{J_{\mu-2,b}f(z)}{J_{\mu-1,b}f(z)},\frac{J_{\mu-3,b}f(z)}{J_{\mu-2,b}f(z)};z\right)
\end{equation*}
is univalent in $\mathbb{U},$ then
\begin{equation}\label{p37}
\Omega \subset\left\{\phi \left(\frac{J_{\mu,b}f(z)}{J_{\mu+1,b}f(z)},\frac{J_{\mu-1,b}f(z)}{J_{\mu,b}f(z)},\frac{J_{\mu-2,b}f(z)}{J_{\mu-1,b}f(z)},\frac{J_{\mu-3,b}f(z)}{J_{\mu-2,b}f(z)};z\right):z \in \mathbb{U}\right\}
\end{equation}
then
\begin{equation*}
 q(z) \prec \frac{J_{\mu,b}f(z)}{J_{\mu+1,b}f(z)}\qquad (z \in \mathbb{U}).
\end{equation*}
\begin{proof}
Let the function $p(z)$ be defined by (\ref{e32}) and $\psi$ by  (\ref{e38}). Since $\phi \in \Phi_{J,2}^{\prime}[\Omega,q],$ (\ref{e39}) and (\ref{p37}) yield 
\begin{equation*}
\Omega \subset \{\psi \left(p(z),zp^{\prime}(z),z^{2}p^{\prime\prime}(z),z^{3}p^{\prime\prime\prime}(z);z\right)z \in \mathbb{U}\}.
\end{equation*}
From equations (\ref{e36}) and (\ref{e37}), we see that the admissible condition for $\phi \in \Phi_{J,2}^{\prime}[\Omega,q]$ in Definition \ref{d9} is equivalent to the admissible condition for $\psi$ as given in Definition \ref{d6} with $ n=2.$ Hence $\psi \in \Psi_{2}^{\prime}[\Omega,q],$ and by using (\ref{p36}) and Lemma \ref{t5}, we have 
\begin{equation*}
q(z)\prec \frac{J_{\mu,b}f(z)}{J_{\mu+1,b}f(z)}\qquad(z \in \mathbb{U}).
\end{equation*}
\end{proof}
\end{thm}
\begin{thm}\label{t19}
Let $\phi \in \Phi_{J,2}^{\prime}[h,q].$ If the function $f \in \mathcal{A}$ and $\frac{J_{\mu,b}f(z)}{J_{\mu+1,b}f(z)} \in \mathbb{Q}_{1}$ and $q \in \mathcal{H}_{1}$ with $q^{\prime}(z)\neq 0$ satisfy the following conditions (\ref{p36}) and

\begin{equation*}
\phi \left(\frac{J_{\mu,b}f(z)}{J_{\mu+1,b}f(z)},\frac{J_{\mu-1,b}f(z)}{J_{\mu,b}f(z)},\frac{J_{\mu-2,b}f(z)}{J_{\mu-1,b}f(z)},\frac{J_{\mu-3,b}f(z)}{J_{\mu-2,b}f(z)};z\right)
\end{equation*}
is univalent in $\mathbb{U},$ then
\begin{equation}\label{p39}
h(z)\prec \phi \left(\frac{J_{\mu,b}f(z)}{J_{\mu+1,b}f(z)},\frac{J_{\mu-1,b}f(z)}{J_{\mu,b}f(z)},\frac{J_{\mu-2,b}f(z)}{J_{\mu-1,b}f(z)},\frac{J_{\mu-3,b}f(z)}{J_{\mu-2,b}f(z)};z\right)
\end{equation}
implies that
\begin{equation*}
 q(z) \prec \frac{J_{\mu,b}f(z)}{J_{\mu+1,b}f(z)}\qquad(z \in \mathbb{U}).
\end{equation*}
\end{thm}

\section{\bf Sandwich-Type Results}
Combining Theorems \ref{t7} and \ref{t14}, we obtain the following sandwich-type theorem.
\begin{cor}\label{c9}
Let $h_{1}$ and $q_{1}$ be analytic functions in $\mathbb{U},$ $h_{2}$ be univalent function in $\mathbb{U},$ $q_{2} \in \mathbb{Q}_{0}$ with $q_{1}(0)=q_{2}(0)=0$ and $\phi \in \Phi_{J}[h_{2},q_{2}]\cap \Phi_{J}^{\prime}[h_{1},q_{1}].$ If the function $f \in \mathcal{A}, J_{\mu+1,b}f(z) \in \mathbb{Q}_{0}\cap \mathcal{H}_{0},$ and  
\begin{equation*}
\phi(J_{\mu+1,b}f(z),J_{\mu,b}f(z),J_{\mu-1,b}f(z),J_{\mu-2,b}f(z);z),
\end{equation*}
is univalent in $\mathbb{U},$ and the condition (\ref{p3}) and (\ref{p30}) are satisfied, then
\begin{equation*}
h_{1}(z)\prec \phi(J_{\mu+1,b}f(z),J_{\mu,b}f(z),J_{\mu-1,b}f(z),J_{\mu-2,b}f(z);z) \prec h_{2}(z)
\end{equation*}
implies that
\begin{equation*}
q_{1}(z)\prec J_{\mu+1,b}f(z)\prec q_{2}(z)\qquad(z \in \mathbb{U}).
\end{equation*}
\end{cor}
Combining Theorems \ref{t10} and \ref{t17}, we obtain the following sandwich-type theorem.
\begin{cor}\label{c10}
Let $h_{1}$ and $q_{1}$ be analytic functions in $\mathbb{U},$ $h_{2}$ be univalent function in $\mathbb{U},$ $q_{2} \in \mathbb{Q}_{1}$ with $q_{1}(0)=q_{2}(0)=1$ and $\phi \in \Phi_{J,1}[h_{2},q_{2}]\cap \Phi_{J,1}^{\prime}[h_{1},q_{1}].$ If the function $f \in \mathcal{A},\frac{J_{\mu+1,b}f(z)}{z} \in \mathbb{Q}_{1}\cap \mathcal{H}_{1},$ and  
\begin{equation*}
\phi \left(\frac{J_{\mu+1,b}f(z)}{z},\frac{J_{\mu,b}f(z)}{z},\frac{J_{\mu-1,b}f(z)}{z},\frac{J_{\mu-2,b}f(z)}{z};z \right),
\end{equation*}
is univalent in $\mathbb{U},$ and the condition (\ref{p17}) and (\ref{p34}) are satisfied, then
\begin{equation*}
h_{1}(z)\prec \phi \left(\frac{J_{\mu+1,b}f(z)}{z},\frac{J_{\mu,b}f(z)}{z},\frac{J_{\mu-1,b}f(z)}{z},\frac{J_{\mu-2,b}f(z)}{z};z \right) \prec h_{2}(z)
\end{equation*}
implies that
\begin{equation*}
q_{1}(z)\prec \frac{J_{\mu+1,b}f(z)}{z}\prec q_{2}(z)\qquad(z \in \mathbb{U}).
\end{equation*}
\end{cor}
Combining Theorems \ref{t12} and \ref{t19}, we obtain the following sandwich-type theorem.
\begin{cor}\label{c11}
Let $h_{1}$ and $q_{1}$ be analytic functions in $\mathbb{U},$ $h_{2}$ be univalent functions in $\mathbb{U},$ $q_{2} \in \mathbb{Q}_{1} $ with $q_{1}(0)=q_{2}(0)=1$ and $\phi \in \Phi_{J,2}[h_{2},q_{2}]\cap \Phi^{\prime}_{J,2}[h_{1},q_{1}].$ If the function $f \in \mathcal{A},$ $\frac{J_{\mu,b}f(z)}{J_{\mu+1,b}f(z)} \in \mathbb{Q}_{1}\cap \mathcal{H}_{1},$ and 
\begin{equation*}
\phi \left(\frac{J_{\mu,b}f(z)}{J_{\mu+1,b}f(z)},\frac{J_{\mu-1,b}f(z)}{J_{\mu,b}f(z)},\frac{J_{\mu-2,b}f(z)}{J_{\mu-1,b}f(z)},\frac{J_{\mu-3,b}f(z)}{J_{\mu-2,b}f(z)};z\right),
\end{equation*}
is univalent  in $\mathbb{U},$ and the condition (\ref{e30}) and (\ref{p36}) are satisfied, then
\begin{equation*}\label{p390}
h(z)\prec \phi \left(\frac{J_{\mu,b}f(z)}{J_{\mu+1,b}f(z)},\frac{J_{\mu-1,b}f(z)}{J_{\mu,b}f(z)},\frac{J_{\mu-2,b}f(z)}{J_{\mu-1,b}f(z)},\frac{J_{\mu-3,b}f(z)}{J_{\mu-2,b}f(z)};z\right)\prec h_{2}(z)
\end{equation*}
implies that
\begin{equation*}
q_{1}(z)\prec \frac{J_{\mu,b}f(z)}{J_{\mu+1,b}f(z)} \prec q_{2}(z)\qquad(z \in \mathbb{U}).
\end{equation*}
\end{cor}
\section*{\bf Acknowledgment}
The present investigation of the second author is supported under the INSPIRE fellowship, Department of Science and Technology, New Delhi, Government of India, Sanction letter No. REL1/2016/2.

\end{document}